\documentclass[reqno,11pt]{amsart}

\usepackage{psfig, amsmath, amsfonts, amssymb, amsthm, amscd}

\DeclareMathSymbol{\leqslant}{\mathalpha}{AMSa}{"36} 
\DeclareMathSymbol{\geqslant}{\mathalpha}{AMSa}{"3E} 
\DeclareMathSymbol{\eset}{\mathalpha}{AMSb}{"3F}     
\renewcommand{\leq}{\;\leqslant\;}                   
\renewcommand{\geq}{\;\geqslant\;}                   
\newcommand{\dd}{\,\text{\rm d}}             

\newcommand{\grad}{\nabla} 
\newcommand{\tigrad}{\tilde\nabla}
\newcommand{\la}{\label} 
 
\newcommand{\be}{\begin{equation}} 
\newcommand{\somme}{\sum_{i=1}^d} 
\newcommand{\sommes}{\sum_{i,j=1}^d} 
\newcommand{\Om}{\Omega}

\newcommand{\sommetn}{\frac{1}{|T_N|}\sum_{x\in T_N}} 
\newcommand{\tipsi}{\tilde\psi} 
\newcommand{\envN}{\Omega_N}
\newcommand{\dhn}[1]{{\rm d}_{{\rm Ham},N} \left\{ #1\right\}}

\def\1{\ifmmode {1\hskip -3pt \rm{I}} \else {\hbox {$1\hskip -3pt \rm{I}$}}\fi}
\newcommand{\df}{\stackrel{\Delta}{=}}
\newcommand{\si}{\sigma}

\newtheorem{Le}{Lemma}[section] 
\newtheorem{Pro}{Proposition}[section] 
\newtheorem{Th}{Theorem}[section] 
\newtheorem{Cor}{Corollary}[section] 
 
\newcommand{\cG}{\ensuremath{\mathcal G}} 
\newcommand{\cH}{\ensuremath{\mathcal H}}

\newcommand{\cL}{\ensuremath{\mathcal L}}

\newcommand{\bbE}{{\ensuremath{\mathbb E}} }

\newcommand{\bbL}{{\ensuremath{\mathbb L}} } 
 
\newcommand{\bbN}{{\ensuremath{\mathbb N}} } 
 
\newcommand{\bbP}{{\ensuremath{\mathbb P}} } 
 
\newcommand{\bbR}{{\ensuremath{\mathbb R}} }

\newcommand{\bbZ}{{\ensuremath{\mathbb Z}} }

\begin{document} 
\title[\,]{Finite volume approximation of the effective diffusion matrix: The case of
independent bond disorder}
\author{Pietro Caputo}
\address{
Dip. Matematica Universita di Roma Tre,
L.go S. Murialdo 1, 00146 Roma, Italy}
\email{caputo@mat.uniroma3.it} 
\author{
Dmitry Ioffe}
\address{
Faculty of Industrial Engineering,
Technion, Haifa 3200, Israel}
\email{ieioffe@ie.technion.ac.il}
\thanks{Partly supported by the ISRAEL SCIENCE FOUNDATION founded by The Israel
Academy of Science and Humanities} 
\date{October 2001}
\vskip 0.2in
\setcounter{page}{1}
\begin{abstract}
Consider uniformly elliptic random walk on $\bbZ^d$ with independent jump rates across
nearest neighbour bonds of the lattice. We show that the infinite volume 
effective diffusion matrix  can be almost surely recovered as the limit of finite 
volume periodized effective diffusion matrices.  
\end{abstract}
\maketitle



\section{Introduction}
Consider the following model of random walk in random environment:
Given $c\geq1$, the space of environments
$$
\Om = \{\xi(b)\in\left[ 1/c ,c \right] \,,\;b\in{\bbZ^d}^*\}\, ,
$$ 
is specified by the jump rates across the nearest neighbour bonds
$b\in{\bbZ^d}^*$. We consider the product measure  
$\mu= \mu_0^{\otimes{\bbZ^d}^*}$ on $\Om$ 
obtained from a probability measure $\mu_0$ supported on $[1/c, c]$. 
For every fixed realization $\xi\in\Om$ of the environment we 
denote $X(t,\xi)$ the nearest neighbour continuous time random walk 
on $\bbZ^d$ which jumps according to 
$\xi(b)$ rates.  The symbol $\bbE_x^\xi$ stands 
for expectation over this process when the walk starts at 
$x\in\bbZ^d$. This is the {\em quenched} regime.  

In the {\em annealed} regime the initial 
random environment is sampled from $\mu$. In the uniformly elliptic case we
consider here the infinite volume annealed homogenization result is well
known \cite{K},\cite{Kozlov},\cite{KV},\cite{DeMasietal},
and the effective 
diffusion matrix $D$ of $X(t)=X(t,\xi )$  
is defined by 
\be 
\forall v\in\bbR^d\qquad (Dv,v) = \sommes v_iD^{ij} v_j = 
\lim_{t\to\infty}\frac{1}{t}\bbE_\mu\bbE^\xi_0\left((X(t,\xi),v)\right)^2 .
\la{vdv} 
\end{equation} 
\smallskip 
We shall briefly recall the Kipnis-Varadhan theory and the variational
formula for $(Dv,v)$ in Subsection~\ref{sub_variational}. Meanwhile let
us proceed with describing the finite volume approximation for the random
environment. This is done in a straightforward fashion: 
Let $N\in\bbZ_+$ and denote by $\pi_N\xi$ 
the periodized bond configuration 
\be 
(\pi_N\xi)_i(x)= \xi_i(\dot x)\,,\;\;x\in\bbZ^d\,,\;\dot x \in T_N\,\df\, 
\bbZ^d/2N\bbZ^d\,. 
\la{perio} 
\end{equation} 
Here $\dot{x}\in\{-N,\dots,N-1\}^d\,,\;x = \dot{x} + 2Nz\,,\;z\in\bbZ^d$,
and $\pi_N\xi$ is a periodic configuration which coincides with 
$\xi$ on the bonds belonging to the torus $T_N$. 
We consider the process $X(t,\pi_N\xi)$ describing a particle 
moving in $\bbZ^d$ according to the rates $\pi_N\xi$. As we shall recall
in the beginning of Section~\ref{section_periodic}, the finite volume 
homogenization result for $X(t,\pi_N\xi)$ is, actually, a much simpler 
statement than the infinite volume one. In particular, there is essentially
no difference between the quenched and annealed regimes, and 
$\forall\xi\in\Om$ the finite volume effective diffusion matrix 
$D_N (\xi )$ is well defined and given by
$$ 
\forall v\in\bbR^d\qquad (D_N (\xi )v,v) = \sommes v_iD_N^{ij}(\xi) v_j = 
\lim_{t\to\infty}\frac{1}{t}\bbE^{\pi_N\xi}_0\left((X(t,\pi_N\xi),v)\right)^2 .
$$ 
Notice that in general the value of $D_N =D_N (\xi )$ depends on the 
realization of the periodized environment $\pi_N\xi$. 

Our main result here asserts that the finite volume effective diffusivities
are rapidly self-averaging and, furthermore, the sequence  
$\{ D_N (\xi )\}$ converges to the annealed infinite volume diffusivity
 $D$.
\begin{Th}
\label{DNtoD}
There exists a power $\nu =\nu (c,d )>0$, such that:
\begin{equation}
\label{concentration_claim}
\mu\left(\left| D_N^{ij} (\xi )-\bbE_\mu D_N^{ij} (\xi )\right| > N^{-\nu}\right)~
<~{\rm e}^{-N^\nu} ,
\end{equation}
for all $N$ sufficiently large. Furthermore,
\begin{equation}
\label{eq_EDNtoD}
\lim_{N\to\infty} \bbE_\mu D_N (\xi )~=~D
\end{equation}
\end{Th}

Of course, as an immediate consequence of the above theorem we obtain 
that 
$\mu$-almost surely,
\begin{equation}
\label{eq_DNtoD}
\lim_{N\to\infty} D_N (\xi )~=~D .
\end{equation}

Our proof of Theorem~\ref{DNtoD} strongly relies on specific properties
of the environment measure $\mu$: the concentration type result 
\eqref{concentration_claim} relies on the independence of the jump rates 
$\xi (b)$, whereas the proof of \eqref{eq_EDNtoD} makes use 
  of the exchangeability of the 
joint distribution of $\xi (b)$-s. Apparently, using the dual variational
description of the inverse matrices $D_N^{-1}$ and $D^{-1}$, the convergence 
result \eqref{eq_DNtoD} could  be established in the general context of 
ergodic jump rates~\cite{Houman}. The finite volume estimate
 \eqref{concentration_claim} should, of course, require 
more stringent 
assumptions on the mixing properties of the $\xi$-field. 
Also, as we shall explain in Subsection~\ref{sub_states}, our proof 
has an additional advantage of giving a natural interpretation 
of the translation covariant states for a harmonic
interface in a random environment in terms of the Funaki-Spohn states 
\cite{FS}.    

We would like to mention that approximation results  of the type 
\eqref{eq_EDNtoD} or \eqref{eq_DNtoD} go beyond 
the general spectral analysis of
the Kipnis-Varadhan approach as developed in \cite{KV} or \cite{DeMasietal}.
 This issue is briefly addressed in 
 Subsection~\ref{sub_spectral}. 
Similar approximation results have been recently derived in \cite{LOV} 
in a much more delicate case 
 of the self-diffusion coefficient of the tagged particle in the exclusion
process. 
\smallskip

Our main motivation for this work 
 came from the theory of massless gradient fields
on $\bbZ^d$. These are specified by the formal Hamiltonian
$$
{\cH}(\phi ) =\sum_x\somme U\left(\nabla_i \phi (x)\right) =
\sum_x\somme U\left(\eta_i (x)\right),
$$
where $\nabla_i$ is the discrete 
lattice gradient; $\eta_i (x)= \nabla_i \phi (x) =\phi (x+ e_i)-\phi (x)$. 
In the uniformly
elliptic case, $1/c \leq U^{\prime\prime} (\cdot )\leq c$,
the infinite volume gradient states exist in any dimension $d\geq 1$, and,
as it 
has been established in the paper by Funaki and Spohn \cite{FS}, any 
translation invariant infinite volume gradient state is decomposable into
the convex combination of the extremal states (the so called Funaki-Spohn
states), which can be constructed as limits of finite volume 
periodized measures with appropriately tilted slopes. The very same 
approximation by the periodized states leads \cite{FS} to a 
meaningful definition  of the (slope-dependent) surface tension. The 
latter is known to be strictly convex \cite{DGI},\cite{GOS}. 

On the mesoscopic level, the integral of the  surface tension happens to be 
precisely the large deviation rate function for the appropriately scaled
height field \cite{DGI}, and it has been conjectured in \cite{GOS} that
the Hessian of the surface tension governs the equilibrium fluctuations
 in the corresponding Funaki-Spohn state. In view of the random walk 
representation developed in \cite{DGI}, this conjecture would follow from
an approximation result along the lines of Theorem~\ref{DNtoD}, but
 for a different type of environment. We refer to \cite{GOS} for further
details. To the best of our knowledge this is still an open problem.
\smallskip

The paper is organized as follows: 
In Section~\ref{section_infinite}
we briefly recall the variational formula for the infinite volume 
effective diffusion matrix and express it in terms of the corrector 
field.   The relevant properties of the corrector field are listed in 
Proposition~\ref{vdvpro} and, other way around, we show that for every
 $v\in\bbR^d$ the corrector field  $\psi_v$ is uniquely determined
by these properties.  Subsequently, $\psi_v$ is recovered as the 
$L^2 (\mu )$ limit of the periodized corrector fields $\tipsi_{N,v}$
 in Proposition~\ref{l2conv} of Section~\ref{section_periodic}. This
already leads to the convergence  of the averaged finite volume effective 
diffusivities $\bbE_\mu  D_N (\xi )$, Corollary~\ref{coro}.  
In Section~\ref{section_as} we show that $D_N (\xi )$ concentrates
around its expected value  $\bbE_\mu  D_N (\xi )$. The concentration 
estimates of Subsection~\ref{sub_concentration} prove
the main claim in Theorem~\ref{DNtoD}.
The key ingredient of the proof is an a-priori $L^p$ estimate of 
Subsection~\ref{sub_lp}, which is used to derive appropriate bounds
on the Hamming distance in Subsection~\ref{sub_hamming}. By 
Meyers argument such $L^p$ estimate follows from discrete version 
of Calderon-Zygmund inequality. Since we could not find direct
 references to the latter, the proof
is sketched in the Appendix.

\section{Diffusivity in the infinite medium} 
\label{section_infinite}
\setcounter{equation}{0}

\subsection{The variational formula}
\label{sub_variational} 
%
%
We will denote $\xi_i(x)$, $i = 1,\dots,d$, $x\in\bbZ^d$, 
the rate $\xi(b)$ at the bond $b = (x,x+e_i)$, with $e_i$ the 
unit vector in the $i$-th direction. 
When we write $\xi_i$ only we mean $\xi_i(0)$. 
$\tau_y\xi$ denotes the shifted configuration 
$(\tau_y\xi)_i(x)=\xi_i(x-y)$, $x,y\in\bbZ^d$, $i = 1,\dots,d$. We also use the 
notation, for any $f: \Om \to \bbR$, 
\[
D_if(\xi) = f(\tau_{-e_i}\xi) - f(\xi)\,,\;\;\;D_i^*f(\xi) = -D_if(\tau_{e_i}\xi) 
\]
In the above notation the associated process of the environment $\xi_t =\tau_{X (t)}\xi$
seen from the particle (see \cite{KV}, \cite{DeMasietal}) is  
a Markov process with generator 
\[ 
\cL f(\xi)= -\somme D_i^*\left(\xi_iD_if\right)(\xi)\,,\;\;\;f:\Om\to\bbR\,. 
\]
By translation invariance $\mu$ is a reversible measure for this process 
and $\cL$ is self adjoint in $L^2(\mu)$.  
\smallskip 
 
It is well known  \cite{DeMasietal}  that   the annealed effective diffusion 
matrix $D$ in \eqref{vdv}  can be recovered from the variational 
formula: 
\be 
(Dv,v) = 2\inf_{f\in L^2(\mu)} 
\somme \bbE_\mu\left(\xi_i(v_i+D_if)^2\right) 
\la{vdvar} 
\end{equation}

\subsection{Corrector field} 
In general, the above variational problem cannot be solved for $f$. However the 
next proposition shows that if we only look at ``gradients'' of $f$ the 
problem has a unique solution. 
\begin{Pro}\la{vdvpro} 
For every $v\in\bbR^d$ we have 
\be 
(Dv,v) = 2\somme v_i \bbE_\mu \xi_i\left(v_i 
+\psi^i_v  \right)
= 2\somme\bbE_\mu\xi_i\left(v_i 
+\psi^i_v  \right)^2\la{vdv3} 
\end{equation} 
where $\psi_v^i$, $i=1,\dots,d$ are the unique elements of $L^2(\mu)$ such that 
\begin{itemize} 
\item{$\bbE_\mu\psi^i_v = 0$\,, 
} 
\item{$\somme D^*_i \left[\xi_i 
\left(v_i + \psi^i_v(\xi)\right)\right] = 0 $\,,\;\;\; 
$\mu$- a.s.\,, 
} 
\item{$D_k\psi^i_v = D_i\psi^k_v$\,,\;\;\; 
$\mu$- a.s.\,,\;\;\;\;$i,k = 1,\dots,d$\,. 
} 
\end{itemize} 
\end{Pro} 
\proof 
For any vector $v\in\bbR^d$ we define the local drift or current along $v$ as
\be 
\varphi_v(\xi) = 
\sum_{i = 1}^d v_i \left(\xi_i(0) - \xi_i(-e_i)\right) 
= - \sum_{i = 1}^d v_iD_i^*\xi_i(0) 
\la{current} 
\end{equation} 
Let us now define the functions 
\be 
\chi^\lambda_v = (-\cL + \lambda)^{-1}\varphi_v\,,\;\;\;\lambda> 0\,,\;v\in\bbR^d\,, 
\la{chilav} 
\end{equation} 
or, in terms of the process: 
\[
\chi^\lambda_v(\xi) = \bbE^\xi_0\int_0^\infty {\rm e}^{-\lambda t} 
\varphi_v(\tau_{X(t)}\xi) \dd t\,. 
\]
Following Kipnis and Varadhan \cite{KV}, Theorem 1.3 - see also 
\cite{Ollapoly} - one can use the spectral resolution of the 
non-negative operator $-\cL$ on $L^2(\mu)$  
to prove 
\be 
\lambda\bbE_\mu(\chi^\lambda_v)^2 \to 0\,,\;\;\;\lambda\to 0\,. 
\la{lachila} 
\end{equation} 
Moreover one can prove that there exist 
functions $\psi^i_v\in L^2(\mu)$, $i = 1,\dots,d$ such that 
\be 
\somme \bbE_\mu\xi_i(D_i\chi^\lambda_v - \psi^i_v)^2 
\to 0\,,\;\;\;\lambda\to 0\,, 
\la{chilatopsi} 
\end{equation} 
which in turn by the ellipticity of $\xi$ implies strong convergence of 
$D_i\chi^\lambda_v$ to $\psi^i_v$ in $L^2(\mu)$. 
In particular, we have 
\[
\bbE_\mu\left(\varphi_v(-\cL)^{-1}\varphi_v\right) 
= \lim_{\lambda\to 0} \bbE_\mu\varphi_v\chi^\lambda_v 
= -\lim_{\lambda\to 0} \somme v_i\bbE_\mu\xi_iD_i\chi^\lambda_v 
= - \somme v_i \bbE_\mu\xi_i\psi^i_v\,. 
\]
Since \eqref{vdvar} is equivalent to
\be 
(Dv,v) = 2\somme \bbE_\mu(\xi_i)v_i^2 
- 2\bbE_\mu\left(\varphi_v(-\cL)^{-1}\varphi_v\right) ,
\la{vdv2} 
\end{equation} 
we obtain the first identity in (\ref{vdv3}).  

To prove the second one we observe that 
\be 
\somme \bbE_\mu\xi_i(\psi^i_v)^2 
= - \somme v_i \bbE_\mu\xi_i\psi^i_v\,. 
\la{vdv5} 
\end{equation} 
This follows from (\ref{lachila}) and 
the convergence of $D_i\chi^\lambda_v$ to $\psi^i_v$ 
(\ref{chilatopsi}) since 
\begin{align*} 
\somme \bbE_\mu\xi_i(\psi^i_v)^2 &= \lim_{\lambda\to 0} 
\somme \bbE_\mu\xi_i(D_i\chi^\lambda_v)^2 = \lim_{\lambda\to 0} 
\bbE_\mu\left(\chi^\lambda_v(-\cL)\chi^\lambda_v\right) 
\nonumber\\ 
&= \lim_{\lambda\to 0} 
\left\{\bbE_\mu \chi_v^\lambda\varphi_v - \lambda\bbE_\mu(\chi^\lambda_v)^2 \right\} 
= - \somme v_i \bbE_\mu\xi_i\psi^i_v\,. 
\end{align*} 
 
Next we show that the functions $\psi_v$ satisfy the three properties in 
the statement. The first property is obvious since $\psi^i_v$ is the $L^2(\mu)$-limit 
of a gradient and $\mu$ is translation invariant. To prove the second 
identity recall that, in view of (\ref{lachila}) 
\[
\lim_{\lambda\to 0}\bbE_\mu(f(-\cL)\chi^\lambda_v) = \bbE_\mu(f\varphi_v)\,,\;\;\;f\in L^2(\mu)\,. 
\] 
Therefore for any $f\in L^2(\mu)$ we have 
\begin{equation*} 
\somme \bbE_\mu fD^*_i \left[\xi_i 
\left(v_i + \psi^i_v(\xi)\right)\right] = 
- \bbE_\mu(f\varphi_v) + 
\lim_{\lambda\to 0} \bbE_\mu(f(-\cL)\chi^\lambda_v) = 0\,. 
\la{laseconda} 
\end{equation*} 
The third identity follows in a similar fashion 
from $D_kD_i\chi^\lambda_v = 
D_iD_k\chi^\lambda_v$, $\lambda>0$. 
 
Finally, to prove uniqueness 
we adapt the argument of Theorem 2 in \cite{PV} and of Proposition 3 in \cite{BD}:  
In the language of \cite{BD} the field  $\psi_v(\xi) =(\psi^1_v(\xi) ,\dots ,
\psi^d_v(\xi))$ is, by the third of the properties in 
\eqref{vdv3},  an $L^2 (\mu )$-cocycle. A straightforward modification of
the proof of Proposition 3 in \cite{BD} reveals that any $L^2 (\mu )$-cocycle 
$(u^1 (\xi ) ,\dots ,u^d (\xi ))$ is, in fact, an $L^2$-limit of gradients:
There exists a sequence 
$g_\gamma\in L^2(\mu)$, $\gamma > 0$, with 
\be 
\bbE_\mu(D_kg_\gamma - u^k)^2 \to 0 
\,,\;\;\gamma\to 0\,,\;\;k=1,\dots,d\,. 
\la{ggamma1} 
\end{equation}
Uniqueness 
of $\psi_v$ is an easy consequence.  Indeed, fix $v$ and let $u^k = \psi^k_{1,v} - 
\psi^k_{2,v}$, with $\psi_{1,v},\psi_{2,v}$ both satisfying the required 
conditions. We then have 
\be 
\sum_{k=1}^d \bbE_\mu\xi_ku^kD_kf = 0\,, 
\;\;\;f\in L^2(\mu) \,. 
\la{ggamma3} 
\end{equation} 
Choosing 
$f = g_\gamma$ in (\ref{ggamma3}), with $g_\gamma$ 
as in (\ref{ggamma1}), we obtain 
\[
\sum_{k=1}^d \bbE_\mu\xi_k(u^k)^2 = 
\lim_{\gamma\to 0} \sum_{k=1}^d \bbE_\mu\xi_ku^kD_kg_\gamma = 0\,, 
\]
which, by the ellipticity of $\xi$, implies $u^k = 0$ $\mu$-a.e., 
$k=1,\dots,d$. \qed

\bigskip 

\noindent 
{\em Remark 1}. 
Using the linearity in $v$ implicit in the definition (\ref{chilav}) 
we may define the functions $\psi^{ij} = \psi^i_{e_j}$, 
$i,j = 1,\dots,d$, so that 
\be 
\psi^i_v = (\Psi v)_i \df \sum_{j=1}^d\psi^{ij}v_j\,. 
\la{psiij} 
\end{equation} 
The diffusivity $D$ is then given by the matrix 
\be 
D = 2\bbE_\mu \Lambda\left(1+\Psi\right) 
\la{Dmatrix} 
\end{equation} 
with $\Lambda$ denoting the diagonal matrix $\Lambda^{ij}=\xi_i\delta_{ij}$. 
The functions $\psi^{ij}$ are often called the corrector fields 
and $\Psi$ the stream matrix.

\smallskip 

\noindent 
{\em Remark 2}. The diffusion constant $D$ has an explicit expression 
if $d=1$. Indeed, for any $v\in\bbR$, $\xi(v+\psi_v(\xi))$ must 
be constant $\mu$-a.s. and normalizing we have 
\be 
\psi_v(\xi) = v\left[\xi^{-1}\left(\bbE_\mu\xi^{-1}\right)^{-1} - 1\right]\,,\;\;\; 
D = \left(\bbE_\mu\xi^{-1}\right)^{-1}\,. 
\la{d=1} 
\end{equation} 
By ellipticity we see that the corrector field $\psi_v$ is uniformly bounded in this case. 
The situation is different if $d\geq 2$. There is no longer an expression for $\psi_v$ and 
all we know apriori is that $\psi\in L^{2}(\mu)$ .
 For $d=2$ this has been upgraded in \cite{Boivin}  to $\psi\in L^{2+\epsilon}(\mu)$
for some $\epsilon> 0$, which  (in two dimensions) leads to an almost sure homogenization
result.  In the sequel, c.f. the remark following Theorem~\ref{Lp_thm}, we 
shall establish that  $\psi\in L^{2+\epsilon}(\mu)$ in any dimension 
$d\geq 2$.
 
\section{Diffusivity in the periodic medium} 
\label{section_periodic}
In the case of the periodized environment \eqref{perio}, 
the process of the environment seen from the particle is  equivalent 
to the projected random walk $\dot X(t,\pi_N\xi)$ on $T_N$, i.e. an 
irreducible Markov chain with finite 
state space and unique invariant measure $\rho_N$ defined by 
\be 
\bbE_{\rho_N}f = \sommetn f(x)\,,\;\;\;f:T_N\to\bbR\,. 
\la{rhon} 
\end{equation} 
Then, introducing the periodic function 
$$\varphi^N_v(x) = \varphi_v(\tau_{-x}\pi_N\xi)\,,$$ 
where $\varphi_v$ is defined by (\ref{current}), 
we arrive as in (\ref{vdv2}) at the expression 
\be 
(D_N(\xi)v,v) = 2\somme \bbE_{\rho_N}((\pi_N\xi)_i)v_i^2 
- 2\bbE_{\rho_N}\left(\varphi^N_v(-\cL_N)^{-1}\varphi^N_v\right)\,, 
\la{vdv20} 
\end{equation} 
where the generator $\cL_N$ describes the jumps on $T_N$, i.e. 
for any periodic function $f$ on $\bbZ^d$ 
\[
\cL_N f(x) = (\pi_N\xi)_i(x)\left[f(x+e_i) - f(x)\right] + (\pi_N\xi)_i(x-e_i)\left[f(x-e_i) - f(x)\right]\,. 
\]
Alternatively, identifying periodic functions with functions on $T_N$ we 
may write 
\[ 
\cL_N f(x) = -\somme \tilde\nabla_i^*\left((\pi_N\xi)_i(x)\tilde\nabla_if(x)\right) 
\,,\;\;\;f:T_N\to\bbR\,, 
\] 
where $\tilde\nabla$ are the discrete gradients on the torus $T_N$.  
We then have the variational principle 
\be 
(v,D_N(\xi)v) = 2\inf_{f: T_N\to\bbR} 
\somme \bbE_{\rho_N}\left((\pi_N\xi)_i(v_i+\tilde\nabla_if)^2\right)\,. 
\la{vdvar0} 
\end{equation} 
In this finite-dimensional setting the above problem may be solved 
directly without the approximation procedure outlined in the proof of 
Proposition \ref{vdvpro}. Namely, define 
\be 
\chi_{N,v}(x) = (-\cL_N)^{-1} \varphi^N_v (x)\,,\;\;\; x\in\bbZ^d 
\la{chinv} 
\end{equation} 
and observe that $\chi_{N,v}$ is a bounded periodic function for each $N$. 
Boundedness follows from the exponential mixing properties of the 
process generated by $\cL_N$, since $\bbE_{\rho_N}\varphi_v^N = 0$. 
As in the previous section we then have 
\be 
(v,D_N(\xi)v) = 2\somme v_i \bbE_{\rho_N} (\pi_N\xi)_i\left(v_i 
+\psi^i_{N,v}  \right) = 2\somme\bbE_{\rho_N} (\pi_N\xi)_i\left(v_i 
+\psi^i_{N,v}  \right)^2 
\la{vdv40} 
\end{equation} 
where $\psi_{N,v}^i$, $i=1,\dots,d$ are defined by 
\be 
\psi_{N,v}^i(x) = \chi_{N,v}(x+e_i) - \chi_{N,v}(x) = \tilde\nabla_i\chi_{N,v}(x) 
\la{psinv} 
\end{equation} 
Clearly $\psi_{N,v}^i$ are periodic functions. The dependence on $\xi$ 
becomes explicit in the notation 
$$\psi_{N,v}^i(x,\xi) = \psi_{N,v}^i(x) = \psi_{N,v}^i(0,\tau_{-x}\pi_N\xi)\,.$$

\subsection{Convergence} 
From now on $\psi_{N,v}^i$ will be regarded as a function of $\xi$, by setting 
\begin{equation}
\label{psi_tilda}
\tipsi_{N,v}^i(\xi) = \psi_{N,v}^i(0,\xi)\,,\;\;\;i=1,\dots,d\,. 
\end{equation}
\begin{Pro}\la{l2conv} 
For each $v\in\bbR^d$, $\tipsi_{N,v}\to\psi_v$ in $L^2(\mu)$, as $N\to\infty$, 
where $\psi_v$ is the field defined in Proposition \ref{vdvpro} 
\end{Pro} 
\proof 
Notice that, as in (\ref{vdv5}) 
\be 
\somme \bbE_{\rho_N}(\pi_N\xi)_i(\psi^i_{N,v})^2 
= - \somme v_i \bbE_{\rho_N}(\pi_N\xi)_i\psi^i_{N,v}\,. 
\la{vdv50} 
\end{equation} 
Then by ellipticity and using Schwarz' inequality 
we obtain
\be 
\max_{i} \bbE_{\rho_N}(\psi^i_{N,v})^2 \leq
 c^2|v|^2 . 
\la{l2bound} 
\end{equation} 
Now we observe that since the bond variables $\xi$ 
are {\em exchangeable} under $\mu$, the distribution of 
$\pi_N\xi$ coincides with that of $\tau_x\pi_N\xi$ for any $x\in\bbZ^d$. 
In particular we see that 
\be 
\bbE_\mu(\tipsi_{N,v}^i)^2 = 
\sommetn \bbE_\mu(\psi^i_{N,v}(0,\tau_x\pi_N\xi))^2 
= \bbE_\mu\bbE_{\rho_N}(\psi^i_{N,v}(\cdot,\xi))^2 
\la{interch} 
\end{equation} 
and (\ref{l2bound}) provides a uniform bound on the $L^2(\mu)$-norm of 
$\tipsi_{N,v}^i$, which implies weak convergence along subsequences. 
Next we show that any weak limit point must satisfy the conditions of Proposition \ref{vdvpro}, 
thus establishing 
the weak convergence of $\tipsi_{N,v}$ to $\psi_v$. 
Let $\tipsi_{\infty,v}^i$ be a weak limit of $\tipsi_{N,v}^i$. 
$\bbE_\mu\tipsi_{\infty,v}^i = 0$ is a consequence of $\bbE_\mu\tipsi_{N,v}^i = 0$, 
which in turn follows from exchangeability. To check the second condition 
we show that for any local function $f\in L^2(\Om)$ we have 
\be 
\somme \bbE_\mu\left[D_if\xi_i 
\left(v_i + \tipsi^i_{\infty,v}(\xi)\right)\right] = 0\,. 
\la{toprove} 
\end{equation} 
From the definition of $\psi^i_{N,v}$, writing out (\ref{chinv}) explicitly 
we have 
\[
\somme \tilde\nabla^*_i\left((\pi_N\xi)_i(x)\psi_{N,v}^i(x,\xi)\right) = 
\varphi^N_v(x) \,, 
\;\;\;x\in T_N\,. 
\]
At $x=0$ this yields 
\be 
\somme\left[ (\tau_{e_i}\xi)_i\left(v_i+\tipsi_{N,v}^i(\tau_{e_i}\pi_N\xi) \right) 
- \xi_i\left(v_i+\tipsi_{N,v}^i(\pi_N\xi) \right)\right] = 0\,. 
\la{at0} 
\end{equation} 
If $f\in L^2(\Om)$ is local, then $f(\xi)= f(\pi_N\xi)$ and $f(\tau_{\pm e_j}\xi) 
= f(\tau_{\pm e_j}\pi_N\xi)$ when $N$ is large, and by the exchangeability 
we have 
\[
\int_{\Om}\mu(\dd\xi)f(\pi_N\xi) 
(\tau_{e_i}\xi)_i\tipsi_{N,v}^i(\tau_{e_i}\pi_N\xi) 
= \int_{\Om}\mu(\dd\xi)f(\tau_{-e_i}\pi_N\xi) \xi_i\tipsi_{N,v}^i(\xi) \,. 
\]
Multiplying (\ref{at0}) with $f$ and integrating 
we see that when $N$ is large 
\[ 
\somme \bbE_\mu \left[D_if \xi_i 
\left(v_i+\tipsi_{N,v}^i\right)\right] = 0\,, 
\]
and the claim (\ref{toprove}) follows from weak convergence. 
The last condition in Proposition \ref{vdvpro} is  proven by similar reasoning, 
using also 
$$\tilde\nabla_i\psi_{N,v}^k= 
\tilde\nabla_i\tilde\nabla_k\chi_{N,v} = 
\tilde\nabla_k\tilde\nabla_i\chi_{N,v} = \tilde\nabla_k\psi_{N,v}^i\,.
$$ 
 
 By the orthogonality relations \eqref{vdv5} and \eqref{vdv50} we, 
using again exchangeability as in \eqref{interch}, infer:
$$
\lim_{N\to\infty} \somme\bbE_\mu \xi_i\left(\tipsi_{N,v}^i\right)^2 ~=~
\somme\bbE_\mu \xi_i\left(\psi_v^i\right)^2 .
$$
In view of the uniform ellipticity condition,  $\xi_i \geq 1/c >0$,  this implies that 
 $\tipsi_{N,v}^i$ converge strongly in $L^2(\mu)$. \qed

\smallskip
 
An immediate corollary of the above results is the mean convergence 
of the diffusion constants $\{D_N (\xi )\}$:
\begin{Cor} 
\la{coro} 
For every $v\in\bbR^d$,
$$
\lim_{N\to\infty} \bbE_\mu \left( v, D_N (\xi )v\right)~=~\left( v, D v\right).
$$
\end{Cor} 
\proof 
Follows  from the representation
formulas \eqref{vdv3} and \eqref{vdv40}, exchangeability and the   $L^2 (\mu )$ convergence
of  $ \tipsi_{N,v}^i$ to $\psi_v^i$. \qed

%

\section{Almost sure convergence of $D_N$} 
\label{section_as}
\setcounter{equation}{0}
In this section we show that the periodized finite volume effective 
diffusion matrices $D_N (\xi )$ converge to the infinite volume 
effective diffusivity $D$ almost surely in the environment, that is 
\eqref{eq_DNtoD} holds.
$L^p$ bounds on the gradients $\tilde{\nabla}\chi_{N,v}$ of the solutions
to \eqref{chinv} set up stage for an application of Talagrand's concentration
estimates, which imply that for every $v\in\mathbb{R}^d$:
\begin{equation}
\label{mu_as_2}
\mu-{\rm a.s}\ \ \ 
\lim_{N\to\infty}\left| \left( D_N (\xi )v, v\right)-{\mathbb E}_\mu
\left( D_N (\xi )v, v\right)\right|~=~0 .
\end{equation}
In the sequel we shall use $\chi_{N,v}^\xi$ instead of $\chi_{N,v}$, just
in order to stress the dependence of solutions of \eqref{chinv} on a 
particular realization $\xi\in [1/c ,c]^{T_N} =\envN$ of the periodized environment.  

\subsection{$L^p$ estimates} 
\label{sub_lp}
The equation \eqref{chinv} is uniformly elliptic in $N\in{\mathbb N}$ and 
$\xi\in\envN$. 
\begin{Th}
\label{Lp_thm}
There exists a power $p =p(c)>2$ and a constant $\alpha =\alpha (c)<\infty$, 
such that uniformly in $v\in {\mathbb R}^d$, $N\in{\mathbb N}$ and 
$\xi\in\envN$,
\begin{equation}
\label{Lp_bound}
\|\tilde{\nabla}\chi_{N,v}^\xi\|_{N,p}~=~\left(\sommetn \left|
 \tilde{\nabla}\chi_{N,v}^\xi(x)\right|^p\right)^{\frac1{p}}~\leq ~
\alpha  | v | .
\end{equation}
\end{Th}
\noindent
{\em Remark.} By the exchangeability of the $\xi$-environment, 
(\ref{Lp_bound}) implies
$$
\sup_N \bbE_\mu \left|\tipsi_{N,v} \right|^p ~<~\infty ,
$$
where $\tipsi_{N,v}$  has been defined in \eqref{psi_tilda}.  Since the infinite volume
corrector field could be obtained as  an almost sure limit of $\tipsi_{N_k,v}$ along a 
subsequence $\{ N_k\}$, it follows that $\psi_v\in L^p (\mu )$ for every $v\in\bbR^d$. 
\vskip 0.1cm
 
Because of the uniform ellipticity of $\cL_N$ the inequality \eqref{Lp_bound}
 follows by Meyers argument \cite{Me}
 (see also Step~1 in Subsection~3.3 in \cite{DGI} for the discrete case)
 once the corresponding $L^q$ estimate holds for the Poisson equation
\begin{equation}
\label{poisson}
\somme \tigrad^*_i\tigrad_i u (x)~=~-\somme \tigrad^*_i f_i (x)
\end{equation}
for some $q >2$. Namely, we have to show that there exists 
$\beta =\beta(q)<\infty$, such that for any $N\in\bbN$ and  every vector field 
$\vec{f}=(f_1,\dots ,f_d )$ on $T_N$ solutions $u$ of \eqref{poisson} 
satisfy:
\begin{equation}
\label{CZ}
\|\tigrad u\|_{N,q}^q~ =\sommetn
|\tigrad u|^q ~\leq ~
\beta (q)\| \vec {f}\|_{N,q}^q .
\end{equation}
\eqref{CZ} is a discrete version of Calderon-Zygmund inequality. We could
not find a direct  reference which would cover the case we consider here.
For the convenience
 of the reader, a brief sketch of the proof is given in the appendix.

\subsection{Bounds on the Hamming distance}
\label{sub_hamming}
For $\xi ,\xi^\prime \in\envN$ let us define the Hamming distance
$$
\dhn{\xi ,\xi^\prime} ~=~\sum_{x\in T_N}\sum_{i=1}^d\delta_{\xi_i (x)\neq 
\xi^\prime_i (x)} .
$$
Set $q=p/(p-2)<\infty$, where $p= p(c)>2$ is as in the statement on Theorem~\ref{Lp_thm}.
We claim that uniformly in $v$ and in $\xi ,\xi^\prime\in\envN$,
\begin{equation}
\label{Ham_bound}
\left|
\left( D_N (\xi )v, v\right) -
\left( D_N (\xi^\prime )v, v\right)\right|~\leq ~
c_1 |v |^2\left[
\frac1{|T_N|} \dhn{\xi ,\xi^\prime}\right]^{1/q}.
\end{equation}
Indeed, by \eqref{vdvar0} and \eqref{vdv40},
\begin{equation}
\label{min_DN}
\begin{split}
\frac12
\left( D_N (\xi )v, v\right)~& =~
\min_{f:T_N\mapsto {\mathbb R}}\sommetn \sum_{i=1}^d
\xi_i (x) \left( v_i +\tilde{\nabla}_i f\right)^2 \\
&\ \ =~
\sommetn \sum_{i=1}^d
\xi_i (x) \left( v_i +\tilde{\nabla}_i
 \chi_{N,v}^\xi\right)^2  ,
\end{split}
\end{equation}
for every $v\in {\mathbb R}^d$ and $\xi\in\envN$. Consequently, 
in view 
of the a-priori estimate  \eqref{Lp_bound}, we infer:
\begin{equation*}
\begin{split}
\frac12 
\big|
( D_N (\xi )v, v) &-
( D_N (\xi^\prime )v, v )\big|\\
&\leq ~
\sommetn \sum_{i=1}^d
\left|\xi_i (x)  -\xi^\prime_i (x)\right| 
\left[
\left( v_i +\tilde{\nabla}_i
 \chi_{N,v}^\xi (x)\right)^2  +
\left( v_i +\tilde{\nabla}_i
 \chi_{N,v}^{\xi^\prime} (x)
\right)^2 
\right]\\
&\ \ \ \ \leq ~c_2  |v |^2
\left[ \sommetn \sum_{i=1}^d 
\left| \xi_i (x) -\xi^{\prime}_{i} (x)\right|^{q}
\right]^{1/q} .
\end{split}
\end{equation*}
Obviously, 
$$
\sum_{x\in T_N} \sum_{i=1}^d 
\left| \xi_i (x) -\xi^\prime_i (x)\right|^{q}~\leq ~
c^{q} \dhn {\xi ,\xi^\prime},
$$
uniformly in $\xi ,\xi^\prime\in\envN$, and the bound \eqref{Ham_bound}
follows.

\subsection{Concentration of $D_N (\xi )$}
\label{sub_concentration}
  By \eqref{min_DN};
$$
0~\leq ~ \left( D_N (\xi ) v ,v\right)~\leq ~ c |v|^2,
$$
for every $v\in {\mathbb R}^d$ and for each realization of the environment
$\xi \in\envN$. Thus, given $\epsilon_N >0$ one can find a value 
$\bar{D}_N (v)\in [0,c|v|^2 ]$, such that 
$\mu \left( A_N (v)\right)\geq 2\epsilon_N/( c |v |^2)$, where we define the 
set $ A_N (v)\subset\envN$ as:
$$
 A_N (v)~=~\left\{ \xi\in\envN~:~\left|\left( D_N (\xi )v,v\right)-\bar{D}_N (v)\right|
<\epsilon_N\right\} .
$$
Every $\xi^\prime\in\envN$ such that 
$$
\left|\left( D_N (\xi^\prime )v,v\right)-\bar{D}_N (v)\right|
>2\epsilon_N ,
$$
certainly satisfies
$$
\min_{\xi\in A_N (v)} \left|\left( D_N (\xi )v,v\right) - \left( D_N (\xi^\prime )v,v\right)
\right| ~>~\epsilon_N .
$$
By \eqref{Ham_bound} we arrive to the following bound on such $\xi^\prime$ in terms of the 
Hamming distance:
$$
\dhn{\xi^\prime ,A_N (v )}~=~\min_{\xi\in A_N (v)}\dhn{\xi^\prime ,\xi}~>~
c_4 \left(\frac{\epsilon_N}{|v|^2}\right)^{q}\, N^d .
$$
We are now in position to use the concentration estimate (5.2) in \cite{Ta}:
\begin{equation}
\label{Talagrand}
\begin{split}
\mu &\left(\xi^\prime~:~\dhn{\xi^\prime ,A_N (v)}  > c_4 
\left(\frac{\epsilon_N}{|v|^2}\right)^{q}\, N^d \right)~\\
&\leq ~\frac{1}{\mu (A_N (v))}
{\rm exp}\left\{ -c_5 \frac{\epsilon_N^{2q}}
{|v|^{4q}}\,
N^d\right\} 
~\leq ~\frac{c|v|^2}{2\epsilon_N}
{\rm exp}\left\{ -c_6(|v|) \epsilon_N^{2q}
\, N^d\right\} 
\end{split}
\end{equation}
It remains to choose $\epsilon_N =\sqrt{N^{-(d-\Delta)/q}}$ 
for some $\Delta>0$: 
Since the sequence of random variables $\left( D_N (\xi )v,v\right)$ is uniformly 
bounded, we readily infer from \eqref{Talagrand} that
$$
\left| {\mathbb E}_\mu \left( D_N (\xi )v,v\right) - \bar{D}_N (v)\right| ~< ~3\epsilon_N
$$
for all large enough values of $N$, and  
the assertion \eqref{concentration_claim} of Theorem~\ref{DNtoD}
follows with $\nu =\min\{\Delta , (d-\Delta)/2q\}$.
 


\subsection{Spectral measures of the local drift $\varphi_v$}
\label{sub_spectral}
We would like to stipulate 
the implication of \eqref{eq_DNtoD} for the spectral properties of the local drifts  
$\varphi_v$ 
as compared to the information which could be extracted from the general Kipnis-Varadhan 
approach: 

In the notation of Section~\ref{section_periodic} let us introduce the empirical measures
$$
\mu_N = \mu_N^\xi ~=~\sommetn \delta_{\tau_{-x}\pi_N\xi} .
$$
For every $\xi\in\Omega$ fixed or, equivalently, for each $\pi_N\xi\in\Omega_N$ fixed there
is an obvious correspondence between the spaces $L^2 (T_N ,\rho_N )$ and 
$L^2 (\Omega_N ,\mu_N )$. By \eqref{vdv2} and \eqref{vdv20} the limiting relation 
\eqref{eq_DNtoD} could be written as 
\begin{equation}
\label{lim_1}
\mu-a.s.\qquad \lim_{N\to\infty}\bbE_{\mu_N}\left(\varphi_v (-\cL_N)^{-1}\varphi_v\right)~=
~\bbE_{\mu}\left(\varphi_v (-\cL)^{-1}\varphi_v\right) .
\end{equation}
Fix now $v\in\bbR^d$. Let $\nu$ and, respectively, $\nu_N =\nu_N^\xi$ be 
the spectral measures of 
$\varphi_v$ relative to the operator
 $( -\cL )$ on $L^2 (\Omega ,\mu )$ and, respectively, relative to the operator $(-\cL_N )$
on $ L^2 (\Omega_N ,\mu_N^\xi )$. Both $-\cL $ and all of $-\cL_N$ are self-adjoint and 
bounded on the respective
spaces. Let $K$ be a common upper bound on the spectral radiuses.  In terms of spectral
measures the limit in \eqref{lim_1} reads as 
\begin{equation}
\label{lim_2}
 \mu-a.s.\qquad \lim_{N\to\infty}\int_0^K \frac{\nu_N^\xi ({\rm d}r)}{r} ~=~
 \int_0^K \frac{\nu({\rm d}r)}{r} .
\end{equation}
Kipnis-Varadhan approach is based on the fact that $\int_0^K \frac{\nu({\rm d}r)}{r}$ (or,
respectively,  $\int_0^K \frac{\nu_N^\xi ({\rm d}r)}{r}$ in the case of the periodized 
environment) is bounded above. We claim that our convergence 
 result \eqref{lim_2} is equivalent to 
the $\mu$-a.s.\ uniform
integrability of the family $\left\{ \nu_N^\xi ({\rm d}r)/r\right\}$.  
Such equivalence  would follow if, for example, we are able to show that $\mu$-almost 
surely the sequence $ \left\{ \nu_N^\xi \right\}$  weakly$\,*$ converges to $\nu$. 
The latter is a consequence of 
\begin{Le}
\begin{equation}
\label{n_limit}
\lim_{N\to\infty} \int_0^K {\rm e}^{-nr}\nu_N^\xi ({\rm d}r)~=~ \int_0^K {\rm e}^{-nr}
\nu ({\rm d}r) ,
\end{equation}
$\mu$-a.s for all $n\in\bbN$. 
\end{Le}
\proof
The claim of the lemma essentially follows from the strong law of large numbers: Since the 
local drift $\varphi_v$ is bounded,
\begin{align}
\label{NM_bulk}
\int_0^K &{\rm e}^{-nr} \nu_N^\xi ({\rm d}r) \\
& =
\frac{1}{[2(N-M)]^d}\sum_{x\in \{-N+M,\dots ,N-M-1\}^d}
\varphi_v (\tau_{-x}\xi )
\bbE_x^\xi \varphi_v^N ( \dot{X} (n ,\pi_N\xi ))~+~{\rm O}\left(\frac{M}{N}\right) , \nonumber
\end{align}
where $\dot{X} (n ,\pi_N\xi )$ is the wrapping  round 
the torus $T_N$ of the random walk 
$X (n ,\pi_N\xi )$  moving in the periodic environment $\pi_N\xi$ and $\bbP_x^\xi$ denotes
the law of such random walk with the starting point $\dot{X} (0 ,\pi_N\xi ) =x$. 
Since the jump
 rates of $X (t,\xi )$ are, uniformly in $\xi\in \Omega$, bounded above by $c$,
$$
\max_{x\in T_N}\sup_{\xi\in\Omega} \bbP_x^\xi\left( \max_{0\leq t\leq n}\left|
X (t,\xi ) -x\right| >M\right)~\leq ~ c_5{\rm e}^{-c_6 M^2 /n} .
$$
Consequently,
\begin{equation}
\label{forget_N}
\bbE_x^\xi \varphi_v^N ( \dot{X} (n ,\pi_N\xi )) ~=~
\bbE_x^\xi \varphi_v ( {X} (n ,\xi )) + {\rm O}\left({\rm e}^{-c_6 M^2 /n}\right),
\end{equation}
uniformly in $x\in \{-N+M,\dots, N-M -1\}^d$. Substituting \eqref{forget_N} into 
\eqref{NM_bulk} and  choosing $M=M(N)=\sqrt{N}$, we arrive to the claim of the lemma.
\qed

\subsection{Some remarks on a massless Gaussian field with bond disorder}
\label{sub_states}
An harmonic interface in the 
quenched random
environment $\xi\in\Om$ is described by the formal Hamiltonian
\be
\cH(\phi) = \frac12\sum_{x}\sum_{i=1}^d\xi_i(x)(\grad_i\phi(x))^2\,.
\la{hami}
\end{equation}
We shall see below how the corrector fields introduced in previous
sections can be used to characterize an interesting class of Gibbs measures
for the interaction (\ref{hami}).

\smallskip

In order to localize the interface in $d=1,2$ we may consider the 
field pinned at the origin, i.e.\ we impose $\phi(0)=0$. Since the interaction
(\ref{hami}) is quadratic a full description of the set $\cG^\xi$ of 
infinite-volume Gibbs measures for a given typical realization $\xi\in\Om$
is available (\cite{Geo}, chap.\ 13). Namely, it is well known that the set of
extremal elements of $\cG^\xi$, denoted ${\rm ext}\cG^\xi$, coincides with the
set of Gaussian fields on $\bbR^{\bbZ^d}$ with a $\xi-$harmonic mean vector 
$\chi: \bbZ^d\to\bbR$ satisfying 
\be
\sum_i \grad^*_i(\xi_i(x) \grad_i\chi(x)) = 0\,,\; x\in\bbZ^d\,;
\quad\chi(0) = 0\,,
\la{chis}
\end{equation}  
and covariance $G^\xi(x,y)$ given by the Green function of the random walk
in the $\xi$-environment killed upon hitting the origin.  ${\rm ext}\cG^\xi$
is thus characterized by solutions to (\ref{chis}).

\smallskip

To study the tilted states associated to (\ref{hami}) it is convenient to 
work directly with the gradient field $\eta_i\df\grad_i\phi$. This is no loss
of information in view of the condition $\phi(0)=0$. In particular,
with the correspondence $\si_i=\grad_i\chi$,
(\ref{chis}) is equivalent to 
\be
\sum_i \grad^*_i(\xi_i(x) \si_i(x)) = 0\,,\quad
\grad_k\si_i(x)=\grad_i\si_k(x)\,,\quad i,k=1,\dots,d\,,\;x\in\bbZ^d\,,
\la{gradchis}
\end{equation}  
and ${\rm ext}\cG^\xi$ is characterized by solutions to (\ref{gradchis}), so that $\nu\in{\rm ext}\cG^\xi$ is a Gaussian measure on $(\bbR^d)^{\bbZ^d}$ with 
mean $\si$ satisfying (\ref{gradchis}) and covariance 
\begin{align*}  
C^\xi_{ij}(x,y)&\df {\rm cov}_\nu(\eta_i(x),\eta_j(y)) \\
& = G^\xi(x+e_i,y+e_j) + G^\xi(x,y) - G^\xi(x+e_i,y) - G^\xi(x,y+e_j)\,.
\end{align*} 

\smallskip

\noindent
{\em Translation covariance and tilted states.} 
A random Gibbs measure is a measurable map $\nu:\,\Om \to \{\cG^\xi,\;\xi\in\Om\}$ such that $\nu^\xi\in\cG^\xi$ for every $\xi\in\Om$. 
The map $\nu^\cdot$ is called {\em translation covariant} when 
$\nu^\xi\circ\theta_x=\nu^{\tau_x\xi}$ for every $\xi\in\Om$, $x\in\bbZ^d$,
where $\theta_x$ denotes the action of translation group and $\tau_x$
is the environment shift. We define the {\em tilted states} for our 
interface in random environment as the set of translation covariant 
random Gibbs measures $\nu^\cdot$ such that $\nu^\xi\in{\rm ext}\cG^\xi$.     
It is not difficult now to use Proposition \ref{vdvpro}
to give a full characterization of tilted states.
Indeed, since $C^\xi_{ij}$ is translation covariant,
we only have to characterize the translation covariant 
maps $\xi\to\si^\xi$ with $\si^\xi$ obeying (\ref{gradchis}). 
For every $v\in\bbR^d$ we know that 
$$\si_v^\xi(x) \df v + \psi_v(\tau_{-x}\xi)$$
has the right properties. On the other hand if we additionally require that 
$\si^\xi(0)$ is in $L^2(\Om,\mu)$ then the uniqueness statement in   
Proposition \ref{vdvpro} allows to conclude that $\si_v^\xi$, $v\in\bbR^d$ are
the only maps with these properties. In other words, for
every  $v\in\bbR^d$ the Gaussian measure $\nu_v^\xi$ 
with covariance $C^\xi$ and
mean $\si_v^\xi$ is the unique tilted state such that 
$\bbE_\mu[\nu_v^\xi(\eta(0))]=v$
and $\bbE_\mu[\nu_v^\xi(|\eta(0)|^2)]<\infty$.
The above measures are the analogue 
in our setting of Funaki-Spohn states for translation invariant
massless gradient fields (\cite{FS}).

\smallskip

\noindent
{\em Tilted states on the torus and surface tension.}
To carry this analogy a little further here we 
mimic the construction of \cite{FS} to describe tilted states as
the infinite volume limit of tilted measure on the torus $T_N$.

Define $\Gamma_N$ as the set of $\tilde\eta\in(\bbR^d)^{T_N}$ such that
$\tilde\eta=\tilde\grad\phi$ for some $\phi\in\bbR^{T_N}$,
$\tilde\grad$ being as usual the discrete gradient on the torus. 
Consider the probability measure on $\Gamma_N$ defined by
\be 
\tilde\nu_{N,v}^\xi(\dd\tilde\eta)=
\frac1{Z_{N,v}^\xi}\exp{
\big(-\frac12
\sum_{x\in T_N}\sum_i(\pi_N\xi)_i(x)(\tilde\eta_i(x)+v_i)^2\big)}
\,m_N(\dd\tilde\eta)
\la{FSN}
\end{equation}
where $m_N(\dd\tilde\eta)$ stands for the image of Lebesgue measure under the map $\phi\to\tilde\grad\phi$. 
Notice that by symmetry $\bbE_{\tilde\nu_{N,0}^\xi}(\tilde\eta(x))=0$, 
$x\in T_N$.
Thus when $v=0$, $\tilde\nu_{N,0}^\xi$ is the centered Gaussian measure 
on $\Gamma_N$ with 
covariance $C_{N,ij}^\xi$ given by the gradients of the Green
function $G_N^\xi(x,y)$
of the periodized random walk on $T_N$ 
with killing upon hitting the origin (\cite{DGI}). Now, the linear tilt 
$v$ only affects the mean in (\ref{FSN}) and it remains to compute 
the mean 
$\bbE_{\tilde\nu_{N,v}^\xi}(\tilde\eta(x))$. We may proceed as follows.
$$
\bbE_{\tilde\nu_{N,v}^\xi}(\tilde\eta(0)) = \int_0^1\frac{\dd}{\dd t}
\bbE_{\tilde\nu_{N,tv}^\xi}(\tilde\eta(0))\,\dd t\,.
$$
But
\begin{align*}
\frac{\dd}{\dd t}
\bbE_{\tilde\nu_{N,tv}^\xi}(\tilde\eta_i(0))&=-
{\rm cov}_{\tilde\nu_{N,tv}^\xi}\big(
\tilde\eta_i(0),\sum_{x\in T_N}
\sum_j(\pi_N\xi)_j(x)(\tilde\eta_j(x)+v_j)v_j
\big)\\
& = - \sum_{x\in T_N}
\sum_j(\pi_N\xi)_j(x)v_jC_{N,ij}(0,x)
=  \chi_{N,v}(e_i)-\chi_{N,v}(0) = \psi^i_{N,v}(0,\xi)
\end{align*}
where in the last line we used the identities (\ref{chinv}) and 
(\ref{psinv})\,.

\smallskip

Now for every $v\in\bbR^d$  let $\nu_{N,v}^\xi$ be the probability measure 
induced by $\tilde\nu_{N,v}^\xi$ on the field $\eta\df\tilde\eta+v$. 
The above computation shows that $\nu_{N,v}^\xi$ is the Gaussian
measure with covariance $C_{N,ij}^\xi$ and mean vector
\be
\bbE_{\nu_{N,v}^\xi}(\eta(x)) = v + \psi_{N,v}(x,\xi)\,.
\la{psinmean}
\end{equation}
From the convergence of $C_{N,ij}^\xi$ we then infer, by Proposition 
\ref{l2conv}, that at least along a subsequence $N'\to\infty$, 
$\nu_{N',v}^\xi$ $\mu-$almost surely converges weakly to the $v-$tilted 
state $\nu_v^\xi$.

\smallskip

A last observation is a simple identity 
connecting the diffusion coefficient to surface tension.
This type of relation was discussed
in \cite{FS,GOS} in the context of translation invariant massless fields. 
Recalling (\ref{FSN}) we define the surface tension on $T_N$
under tilt $v\in\bbR^d$
as 
$$
\si_N^\xi(v)= 
\frac1{|T_N|}\log{\frac{Z_{N,0}^\xi}{Z_{N,v}^\xi}}\,.
$$
By performing the Gaussian integrals and using identities
(\ref{chinv}) and (\ref{vdv20}) one obtains the sought relation
\be
\si_N^\xi(v)= \frac14 (D_N(\xi)v,v)\,.
\la{surfacet}
\end{equation}
An immediate corollary of Theorem \ref{DNtoD} is then the a.s.\ convergence
of $\si_N^\xi(v)$ to $\frac14 (Dv,v)$. As we have 
already mentioned, the corresponding problem 
raised in \cite{GOS} in the context of an-harmonic models  remains open.

\section{Appendix}
\setcounter{equation}{0}
The proof of Theorem~\ref{Lp_thm} is a straightforward adjustment to the discrete Poisson 
equation \eqref{poisson} on $T_N$  of the arguments 
employed in Chapter~9 of \cite{GT}: Given $N\in\bbN$ and a vector field
$\vec{g}$ on $T_N$  define the distribution function $\eta_{N,\vec{g}}$
 of $\vec{g}$ as:
$$
\eta_{N,\vec{g}}(t)~=~\frac1{|T_N|}\#\left\{ x\in T_N~:~|\vec{g} (x)|>t\right\} .
$$
By the Marcinkiewicz interpolation theorem applied for the solution map 
$T:\vec{f}\mapsto \tigrad u$ it would be enough to check that there exist
 constants $c_1$ and $c_2$ such that for every $t\geq 0$
\begin{equation}
\label{L_1_bound}
\eta_{N, \tigrad u}(t)~\leq ~ c_1 \frac{\|\vec{f}\|_{N,1}}{t} , 
\end{equation}
and, 
\begin{equation}
\label{L_2_bound}
\eta_{N, \tigrad u}(t)~\leq ~ c_2 \frac{\|\vec{f}\|_{N,2}^2}{t^2} , 
\end{equation}  
Since $\|\tigrad u\|_{N,2}\leq \|\vec{f}\|_{N,2}$, the second estimate
 trivially follows from Markov inequality (with $c_2 =1$). It remains to
 check \eqref{L_1_bound}

To facilitate the exposition we shall consider only the dyadic case $N=2^n$.
Also, it is enough to derive the bound only for the vector fields $\vec{f}$
of the form $\vec{f}(x) = f(x)\vec{e}_1$. 

For every $k=0,1,\dots,n-1$ define the decomposition
$$
T_N~=~\bigvee_{x\in KT_{N/K}} B_x^k
$$
of the torus $T_N$ into the lattice boxes $B_x^k = x +\{0,\dots, 2^{k}-1\}^d$
of the linear size $K=2^k$. 
The boxes on different $k$-scales are naturally ordered by the inclusion: we
say $B_y^{l}$ is a predecessor  of $B_x^k$ if $l>k$ and $B_x^k\subset B_y^l$. By
 definition $T_N$ is a predecessor  of any of the boxes $B_x^k$ on 
any of the $k$-scales.  Let us fix a number $t>\|f\|_{N ,1}$. We say that 
a box $B_x^k$ is {\bf correct} if 
$$
\frac1{|B_x^k |}\sum_{y\in B_x^k} |f (y)|~=~\frac1{K^d}
\sum_{y\in B_x^k} |f (y)| ~\leq ~t.
$$
By the choice of $t$ the common ancestor $T_N$ is always correct. An incorrect
box $B_x^k$ is called {\bf marked} if it is incorrect, but all his predecessors
are correct. By construction,
\begin{equation}
\label{marked_bound}
t~<~ \frac1{|B_x^k |}\sum_{y\in B_x^k} |f (y)|~\leq ~2^dt,
\end{equation}
for every marked box $B_x^l$. Since all the marked boxes are disjoint,
 the first of the inequalities in \eqref{marked_bound} implies that the 
total marked volume
\begin{equation}
\label{marked_volume}
\frac{1}{|T_N|}\sum_{(k,x) :B_x^k\ \text{is marked}}|B_x^k |~\leq ~
\frac{\|f\|_{N,1}}{t} .
\end{equation}
Let us renumber all the marked boxes as $B_1 ,\dots, B_l$ and decompose
$f$ as 
$$
f(x) = f(x)\1_{T_N\setminus\cup B_i} (x)+\sum_{i=1}^lf(x)\1_{B_i}(x)~\df~
f^0 (x) + \sum_{i=1}^lf^i (x) .
$$ 
Accordingly, we decompose solutions $u$ of \eqref{poisson} as 
$u = u^0 +\sum_1^l u^j$, where
\begin{equation}
\label{equation_B_i}
\somme \tigrad^*_i\tigrad_i u^j (x)~=~-\tigrad^*_1 f^j (x),
\end{equation}
for $j=0,1,\dots,l$.  Evidently,
$$
\eta_{N ,\tigrad u}(2t )~\leq ~\eta_{N ,\tigrad u^0}(t ) + 
\eta_{N ,\tigrad \sum_1^lu^i}(t ) .
$$
Since $|f^0 |$ is bounded above by $t$, it follows from \eqref{L_2_bound}
 that 
\begin{equation}
\label{smooth_part}
\eta_{N ,\tigrad u^0}(t )~\leq ~\frac{\|f^0\|_{N,2}^2}{t^2}~\leq ~
\frac{\|f^ 0 \|_{N,1}}{t} .
\end{equation}
The bulk of the work is, thus, to derive the $L^1$  estimate on the
 distribution function $\eta_{N ,\tigrad \sum_1^lu^i}$ corresponding to 
the irregular part $\sum_1^l f^l$ of the vector field $\vec{f} = f\vec{e}_1$.
The equation \eqref{equation_B_i} feels the right hand side $f^i$ only 
inside the box $B_i$. For $B_i = B_x^k$ let us define the enlargement 
$\bar{B}_i =\bar{B}_x^k$ via
$$
\bar{B}_i = \bar{B}_x^k = \bigcup_{y\in KT_{N/K}:\|y-x\|_\infty \leq K}B_y^k .
$$
In other words, $\bar{B}_x^k$ is the union of $B_x^k$ with all the 
$*$-contiguous boxes on the $k$-th scale. By \eqref{marked_volume},
\begin{equation}
\label{marked_volume_bar}
\frac1{|T_N|}\sum_{i=1}^{l} |\bar{B}_i |~\leq ~ 3^d\frac{\| f\|_1}{t} .
\end{equation}
In order to estimate $\eta_{N,\tigrad u^i}$ outside $\bar{B}_i$ write:
\begin{equation}
\label{grad_u_i}
\tigrad_j u^i (x) = \sum_{y\in B_i} \tigrad_1\tigrad_j G_N (x-y ) f^i (y) ~\df ~
\sum_{y\in B_i} \tigrad^2_{1,j}G_N (x-y ) f^i (y)
\end{equation}
where
$$
\tigrad_{1,j}^2 G_N (z ) = \sum_{m=0}^\infty  \tigrad_{1,j}^2 p_m (z),
$$
and $p_m$ is the $m$-step transition function of the simple random walk on 
$T_N$. 

There is no loss to assume that $f^i$ has zero average:
$$
\sum_{y\in B_i} f^i (x)~=~0.
$$
Thus,  for $B_i =B_{y_0}^k$,  we can rewrite \eqref{grad_u_i} as
\begin{equation}
\label{grad_u_i_repr}
\tigrad_j u^i (x)~=~\sum_{y\in B_i} \left( \tigrad^2_{1,j}G_N (x-y ) 
-\tigrad^2_{1,j}G_N (x-y_0 )\right)  f^i (y) .
\end{equation}
By Theorem~1.5.5 in \cite{La} ($d\geq 3$) or by Theorem~1.6.5 in \cite{La} ($d=2$),
$$
\left|\tigrad^2_{1,j}G_N (x-y ) 
-\tigrad^2_{1,j}G_N (x-y_0 )\right| ~\leq ~c_3 (d) \frac{K}{|x-y_0 |^{d+1}} .
$$
It follows that
$$
\sum_{x\in T_N\setminus\bar{B}_i}\left| \tigrad u^i\right| ~\leq ~
c_4\sum_{y\in B_i} |f (y)| ,
$$
for all $i=1,\dots ,l$. Since the marked boxes $B_1 ,\dots ,B_l$ are disjoint,
$$
\sum_{x\in T_N\setminus\cup\bar{B}_i}\left| \tigrad \sum_{i=1}^l u^i\right|~
\leq ~\sum_{y\in T_N} |f (x)| .
$$
As a result,
$$
\frac1{|T_N|}\#\left\{ x\in T_N\setminus\cup\bar{B}_i~:~
\left| \tigrad \sum_{i=1}^l u^i\right| >t\right\}~\leq ~c_4\frac{\|f\|_{N,1}}{t},
$$ 
 which, by \eqref{marked_volume_bar}, leads to the desired estimated on the distribution
function of the $\tigrad \sum_1^l u^i$ part:
$$
\eta_{N ,\tigrad \sum_1^lu^i}(t )~\leq ~c_5 \frac{\|f\|_{N,1}}{t} .
$$
The proof of Theorem~\ref{Lp_thm} is concluded.

\end{document}